\documentclass{amsart}
\usepackage{graphicx}
\usepackage[top=1.5in,bottom=1.5in,left=1.5in,right=1.5in]{geometry}
\newtheorem{theorem}{Theorem}[section]
\newtheorem*{theorem A}{Theorem A}
\newtheorem*{theorem B}{N\"olker's Theorem}
\newtheorem{lemma}{Lemma}[section]
\newtheorem{proposition}{Proposition}[section]
\newtheorem{corollary}{Corollary}[section]

\theoremstyle{remark}

\theoremstyle{remark}

\theoremstyle{definition}

\newtheorem{definition}{Definition}[section]

\numberwithin{equation}{section}
\def\({\left ( }
\def\){\right )}
\def\<{\left < }
\def\>{\right >}

\begin{document}

\title{ CERTAIN RESULTS ON $N(\kappa)$-CONTACT METRIC MANIFOLDS}
\author{ Absos Ali Shaikh$^{1}$ and Sunil Kumar Yadav$^{2}$}
%    Address of record for the research reported here
\address{$^{1}$ Department of Mathematics, University of Burdwan, Burdwan-713 104, West Bengal, India}
\email{aashaikh@math.buruniv.ac.in, aask2003@yahoo.co.in}

\address{$^{2}$ Department of Mathematics, Poornima College of Engineering, Sitapura, Jaipur-302020, Rajasthan, India}
\email{prof\_sky16@yahoo.com}

%    Address of record for the research reported here


%%%%%%%%%%%%%%%%%%%%%
\begin{abstract}
In this paper, $N(\kappa)$-contact metric manifolds satisfying the conditions $\widetilde{C}(\xi,X)\cdot\widetilde{C}=0$, $\widetilde{C}(\xi,X)\cdot R=0$, $\widetilde{C}(\xi,X)\cdot S=0$, $\widetilde{C}(\xi,X)\cdot C=0$, $C\cdot S=0$ and $R\cdot C=f_{C}Q(g,C)$ have been investigated and obtained their classification. Among others it is shown that a Weyl-pseudosymmetric $N(\kappa)$-contact metric manifold is either locally isometric to the Riemannian product $E^{n+1}(0)\times S^{n}(4)$ or an $\eta$-Einstein manifold. Finally, an example is given.
\end{abstract}
\maketitle
\section{Introduction}
\noindent Let $M$ be an $m (\geq 3)$-dimensional connected smooth Riemannian manifold endowed with the Riemannian metric $g$, Levi-Civita connection $\nabla$. Let $R$, $S$, $r$ be respectively the Riemannian curvature tensor, Ricci tensor and scalar curvature of $M$ such that 

\begin{equation}\label{GrindEQ__1_1_} 
R(X,Y)=[\nabla_{X},\nabla_{Y}]-\nabla_{[X,Y]}, 
\end{equation}

\begin{equation}\label{GrindEQ__1_2_} 
(X\wedge_g Y)Z=g(Y,Z)X-g(X,Z)Y, 
\end{equation}

\begin{equation}\label{GrindEQ__1_3_} 
g(QX,Y)=S(X,Y), 
\end{equation}

\begin{equation}\label{GrindEQ__1_4_} 
S^{2}(X,Y)=S(QX,Y). 
\end{equation}

\noindent for any vector field $X$, $Y$, $Z$ $\in \chi (M)$, $\chi (M)$ being the Lie algebra of all smooth vector fields on $M$. Then the concircular curvature $ \widetilde{C}$ and Weyl-conformal curvature tensor $C$ on $M$ are given by \cite{YAN}

\begin{equation}\label{GrindEQ__1_5_} 
\widetilde{C}(X,Y)=\(R-\frac{r}{m(m-1)}\wedge_g\)\left(X, Y\right),
\end{equation}

\begin{equation}\label{GrindEQ__1_6_} 
C(X,Y)=R(X,Y)-\frac{1}{m-2}\left\{ (X\wedge QY)+(QX \wedge Y)-\frac{r}{m-1}(X\wedge Y)\right\}.
\end{equation}
\noindent For a $(0,l)$-tensor field $T$, $l\geq1$, on $(M,g)$, the tensors $R\cdot T$ and $Q(g,T)$ are defined as

 \begin{align}
\label{GrindEQ__1_7_}
 (R(X,Y) \cdot T) (X_{1}, \cdots , X_{l}) = & -T(R(X,Y) X_{1},X_{2},  \cdots ,X_{l})-\cdots \\ & \cdots -T(X_{1}, \cdots ,X_{l-1}, R(X,Y) X_{l}),\nonumber
 \end{align}

\begin{align} \label{GrindEQ__1_8_}
Q(g,T)(X_1, \cdots ,X_l;X,Y)= & -T((X\wedge Y)X_1,X_2, \cdots ,X_l)-
\cdots \\ & \cdots -T(X_1, \cdots ,X_{l-1},(X\wedge Y)X_l), \nonumber
\end{align}

\noindent respectively \cite{DES2}.

\noindent If the tensors  $R\cdot C$  and  $Q(g,C)$  are linearly dependent then $M$ is called Weyl-pseudosymmetric, that is,
\begin{equation}\label{GrindEQ__1_9_} 
 R\cdot C=f_CQ(g,C),
\end{equation}
holding on the set $U_C=\left\{x\in M\left|C\right.\neq 0 \,\,at\,\, x\right\}$, where $f_C$ is some function on $U_C$. If $R\cdot C=0$ then $M$ is called Weyl-semisymmetric (\cite{DES1}, \cite{DES2}, \cite{DES3}). If $\nabla C=0$ then $M$ is called conformally symmetric \cite{CHA}. It is well-known  that a conformally symmetric manifold is Weyl-semisymmetric. Furthermore, we define the tensor $C\cdot S$ on $M$
\begin{equation}\label{GrindEQ__1_10_} 
(C(X,Y)\cdot S)(Z,W)=-S(C(X,Y)Z,W)-S(Z,C(X,Y)W).
\end{equation}

\noindent A contact metric manifold $M$ of dimension $m(=2n+1), (n>1)$ is a quadruple $(\phi, \xi, \eta, g)$ such that 
\begin{equation}\label{GrindEQ__1_11_} 
\eta(X)=g(X,\xi),\,\,\,\,d\eta(X,Y)=g(X,\phi Y),\,\,\,\,\phi^{2}(X)=-X+\eta(X)\xi, 
\end{equation}
\begin{equation}\label{GrindEQ__1_12_} 
\phi\xi=0,\,\,\,\,\eta\circ\phi=0,\,\,\,\, g(\phi X,\phi Y)=g(X,Y)-\eta(X)\eta(Y), 
\end{equation} 

\noindent  for all $X,Y \in \chi(M)$.

\noindent On a contact metric manifold $M$, a tensor field $h$ is defined by $h=\frac{1}{2} \ell _{\xi } \phi $, 
where $\ell$ denotes the operator of Lie differentiation. Then $h$ is symmetric and satisfies (\cite{BLA1})
\begin{equation}\label{GrindEQ__1_13_} 
h\xi=0,\,\,\,h\phi=-\phi h,\,\,\,Tr.h=Tr.\phi h=0. 
\end{equation} 
\begin{equation}\label{GrindEQ__1_14_} 
\nabla_{X}\xi=-\phi X-\phi hX.
\end{equation}

\noindent If $\xi$ is a Killing vector field then $M$ is said to be a $K$-contact metric manifold, and $M$ is Sasakian if and only if
\begin{equation}\label{GrindEQ__1_15_} 
R(X,Y)\xi=\eta(Y)X-\eta(X)Y.
\end{equation}

\noindent The notion of $\kappa$-nullity distribution was introduced by Tanno \cite{TAN} for a real number $\kappa$ as a distribution 
\begin{equation} \label{GrindEQ__2_1_}
N(\kappa):p\rightarrow N_{p}(\kappa)=\left\{Z \in T_{p}M:R(X,Y)Z=\kappa[g(Y,Z)X-g(X,Z)Y]\right\},
\end{equation}
\noindent for any $X,Y \in T_{p}M $. Hence if $\xi \in N(\kappa)$ then 
\begin{equation} \label{GrindEQ__2_2_}
R(X,Y)\xi=\kappa[\eta(Y)X-\eta(X)Y]
\end{equation}
 holds. Thus a contact metric manifold $M$ for which $\xi \in N(\kappa)$ is called a $N(\kappa)$-contact metric manifold. From (1.15) and (1.17) it follows that a $N(\kappa)$-contact metric manifold is a Sasakian if and only if $\kappa=1$. On the other-hand if $\kappa=0$, then the manifold is locally isometric to the product $E^{n+1}(0)\times S^{n}(4)$ for $n>1$ and flat for $n=1$ \cite{BLA2}. Also in a $N(\kappa)$-contact metric manifold, $\kappa$ is always a constant such that $\kappa \leq1$ \cite{TAN}. Throughout the paper by $M$ we mean a $(2n+1),(n>1)$-dimensional $N(\kappa)$-contact metric manifold unless otherwise stated. 
 
 The paper is structured as follows: Section 2 deals with some requisitory curvature properties of $N(\kappa)$-contact metric manifold. Section 3 is concerned with main results and it is shown that a $(2n+1)$-dimensional $N(\kappa)$-contact metric manifold $M$ satisfies $\widetilde{C}(\xi,X)\cdot\widetilde{C}=0$  (resp., $\widetilde{C}(\xi,X)\cdot C=0$, $\widetilde{C}(\xi,X)\cdot S=0$) if and only if the manifold is either $N(1-\frac{1}{n})$-contact metric manifold or locally isometric to the hyperbolic space $H^{2n+1}(-\kappa)$ (resp., Einstein manifold, $\eta$-Einstein manifold). Also it is shown that if $M$ is a Weyl-pseudosymmetric then it is either locally isometric to the Riemannian product $E^{n+1}(0)\times S^{n}(4)$ for $n>1$ and flat for $n=1$ or $\eta$-Einstein manifold.

\begin{center}
\section { $N(k)$-Contact Metric Manifolds}
\end{center}

Generalizing the notion of $N(\kappa)$-contact metric manifold in 1995, Blair, Koufogiorgos and Papantoniou \cite{BLA3} introduced the notion of $N(\kappa,\mu)$-contact metric manifold, for real numbers $\kappa$ and $\mu$ as a distribution :
\begin{equation} \label{GrindEQ__2_3_}
N(\kappa,\mu):p\rightarrow N_{p}(\kappa,\mu)=\left\{Z \in T_{p}M:R(X,Y)Z=\kappa[g(Y,Z)X-g(X,Z)Y]\right.
\end{equation} 
 $${\,\,\,\,\,\,\,\,\,\,\,\,\,\,\,\,\,\,}+ \mu[g(Y,Z)hX-g(X,Z)hY],$$

\noindent for any $X,Y \in T_{p}M$. Hence if the characteristic vector field $\xi$ belongs to the $(\kappa,\mu)$-nullity distribution, then 
\begin{equation} \label{GrindEQ__2_4_}
R(X,Y)\xi=k [\eta(Y)X-\eta(X)Y]+\mu [\eta(Y)hX-\eta(X)hY].
\end{equation}

\noindent A contact metric manifold $M$ satisfying the relation (2.2) is called a $N(\kappa,\mu)$-contact metric manifold or simply a $(\kappa,\mu)$-contact metric manifold. In particular, if $\mu = 0$, then the relation (2.2) reduces to (1.17) and hence a $N(\kappa)$-contact metric manifold is a $N(\kappa, 0)$-contact metric manifold.

\noindent Let $M$ be a $N(\kappa)$-contact metric manifold. Then the following relations hold (\cite{SHA}, \cite{TAN}, \cite{SHA4}, \cite{SHA6}, \cite{SHA7}):
\begin{equation} \label{GrindEQ__2_5_}
Q \phi-\phi Q=4(n-1)h\phi,
\end{equation}
\begin{equation} \label{GrindEQ__2_6_}
h^{2}=(\kappa-1)\phi^{2},\,\,\kappa \leq1,
\end{equation}
\begin{equation} \label{GrindEQ__2_7_}
Q\xi=2n\kappa \xi,
\end{equation}
\begin{equation} \label{GrindEQ__2__8_}
R(\xi,X,Y)=\kappa [g(X,Y)\xi-\eta(Y)X].
\end{equation}
\noindent In view of (1.1) and (1.2), it follows from (2.1), (2.2), (2.3) and (2.4) that in a $N(\kappa)$-contact metric manifold, the following relations hold:
\begin{equation} \label{GrindEQ__2_9_}
Tr.h^{2}=2n(1-\kappa),
\end{equation}
\begin{equation} \label{GrindEQ__2_10_}
S(X,\phi Y)+S(\phi X,Y)=2(2n-2)g(\phi X,hY),
\end{equation}
\begin{equation} \label{GrindEQ__2_11_} 
S(\phi X,\phi Y)=S(X,Y)-2nk\eta(X)\eta(Y)-2(2n-2)g(hX,Y),
\end{equation} 
\begin{equation} \label{GrindEQ__2_12_} 
S(X,\xi)=2n\kappa \eta(X),
\end{equation} 
\begin{equation} \label{GrindEQ__2_13_} 
Q\phi+\phi Q=2\phi Q+2(2n-2)h\phi,
\end{equation} 
\begin{equation} \label{GrindEQ__2_14_} 
\eta(R(X,Y)Z)=\kappa [g(Y,Z)\eta(X)-g(X,Z)\eta(Y)],
\end{equation}
\begin{equation} \label{GrindEQ__2_15_} 
S(\phi X,\xi)=0,
\end{equation}

\noindent for any vector field $X,Y\in \chi (M)$. Also in a $N(\kappa)$-contact metric manifold the scalar curvature $r$ is given by (\cite{BLA3}, \cite{SHA})
\begin{equation} \label{GrindEQ__2_16_}
r=2n(2n-2+\kappa).
\end{equation}

\noindent In \cite{BOE} Boeckx introduced an invariant on a non-Sasakian contact metric manifold as 
\begin{equation*} \label{GrindEQ__2_17_}
I_{M}=\frac{1-\frac{\mu}{2}}{\sqrt{1-\kappa}}
\end{equation*}
and showed that for two non-Sasakian $(\kappa,\mu)$-manifolds $M_{1}$ and $M_{2}$, we have $I_{M_{1}}=I_{M_{2}}$ if and only if up to a $D$-homothetic deformation, the two manifolds are locally isometric as contact metric manifolds.
Thus, we see that from all non-Sasakian $(\kappa,\mu)$-manifolds of dimension
$(2n + 1)$ and for every possible value of the invariant $I_{M}$, one $(\kappa,\mu)$-manifold $M$ can be obtained. For $I_{M}>-1$ such examples may be found from the standard contact metric structure on the tangent sphere bundle of a manifold
of constant curvature $c$ where we have $I_{M}=\frac{1+c}{1-c}$. Boeckx also gives a Lie algebraic construction for any odd dimension and value of $I_{M}<-1$.

 Using this invariant, Blair, Kim and Tripathi \cite{BLA4} constructed an example of a $(2n + 1)$-dimensional $N(1-\frac{1}{n})$-contact metric manifold $n>1$. Since the Boeckx invariant for a $(1-\frac{1}{n},0)$-manifold is $\sqrt{n}>-1$, we consider the tangent sphere bundle of an $(n+1)$-dimensional manifold of constant curvature $c$ so chosen that the resulting $D$-homothetic deformation will be a $(1-\frac{1}{n},0)$-manifold. That is, for $k=c(2-c)$ and $\mu=-2c$, we solve
\begin{equation*} \label{GrindEQ__2_18_}
1-\frac{1}{n}=\frac{\alpha +a^{2}-1}{a^{2}},\,\,\,\,0=\frac{\mu+2a-2}{a},
\end{equation*}
\noindent for $a$ and $c$. We have
\begin{equation*} \label{GrindEQ__2_19_}
c=\frac{\sqrt{n}\pm1}{n-1},\,\,\,\,a=1+c,
\end{equation*}
and taking $c$ and $a$ to be these values we obtain $N(1-\frac{1}{n})$-contact metric manifold.

\begin{definition} A $N(\kappa)$-contact metric manifold $M$ is said to be $\eta$-Einstein if its Ricci tensor $S$ of type $(0,2)$ is of the form
\begin{equation} \label{GrindEQ__2_20_} 
S(X,Y)=c_{1}g(X,Y)+c_{2}\eta (X)\eta(Y),
\end{equation} 
\end{definition}
\noindent where $c_{1},c_{2}$ are smooth functions on $M$. Contracting (2.15), we have
\begin{equation} \label{GrindEQ__2_21_} 
r=(2n+1)c_1+c_2.
\end{equation} 
Beside this, taking $X=Y=\xi$ in (2.15) and using (2.10) we also have
\begin{equation} \label{GrindEQ__2_22_} 
2n\kappa=c_1+c_2.
\end{equation} 
Hence in view of (2.15),(2.16) and (2.17), we have the result.
\begin{proposition}\label{pro4.1}
In an $\eta$-Einstein $N(\kappa)$-contact metric manifold $M$, the Ricci tensor $S$  is of the form
\begin{equation} \label{GrindEQ__2_23_} 
S(X,Y)=\left(\frac{r}{2n}-\kappa\right)g(X,Y)+\left((2n+1)\kappa-\frac{r}{2n}\right)\eta(X)\eta(Y).
\end{equation} \end{proposition}
\begin{proposition}\label{pro3.1}
\cite{BLA2} A contact metric manifold $M$ satisfying the condition $R(X,Y)\xi=0$ for all $X,Y$ is locally isometric to the Riemannian product of a fat $(n+1)$-dimensional manifold and an $n$-dimensional manifold of positive curvature $4$, i.e., $E^{n+1}(0)\times S^{n}(4)$ for $n>1$ and flat for $n=1$.
\end{proposition}

\begin{proposition}\label{pro3.1}
\cite{BLA6} Let $M$ be an $\eta$-Einstein manifold of dimension $(2n+1)$, $(n\geq1)$, if $\xi$ belongs to the $\kappa$-nullity distribution, then $\kappa=1$ and the structure is Sasakian.
\end{proposition}

\begin{proposition}\label{pro3.1}
\cite{TAN} Let $M$ be an Einstein manifold of dimension $(2n+1)$, $(n\geq2)$, if $\xi$ belongs to the $\kappa$-nullity distribution, then $\kappa=1$ and the structure is Sasakian.
\end{proposition}

\begin{center}
\section{Main results}
\end{center}
In this section, we focus on the characterization of $N(\kappa)$-contact metric manifolds satisfying the condition $\widetilde{C}(\xi,X)\cdot\widetilde{C}=0$, $\widetilde{C}(\xi,X)\cdot R=0$, $\widetilde{C}(\xi,X)\cdot S=0$, $\widetilde{C}(\xi,X)\cdot C=0$, $C\cdot S=0$ and $R\cdot C=f_{C}Q(g,C)$ and deduce some results. For equivalency of several curvature restrictions on a semi-Riemannian manifold, we refer to the reader to see \cite{SHA1}, \cite{SHA2}.
\begin{theorem}
\label{thm4.1}
A $(2n+1)$-dimensional $N(\kappa)$-contact metric manifold $M$, satisfies
\begin{equation*} \label{GrindEQ__3_1_} 
\widetilde{C}(\xi,X)\cdot\widetilde{C}=0
\end{equation*}
if and only if either the manifold is $N(1-\frac{1}{n})$-contact metric manifold or it is locally isometric to the hyperbolic space $H^{2n+1}(-\kappa)$.
\end{theorem}

\noindent\textbf{Proof.} In view of (1.17) and (2.6), equation (1.5) reduces
\begin{equation} \label{GrindEQ__3_2_} 
\widetilde{C}(\xi,Y)Z=\left\{\kappa-\frac{r}{2n(2n+1)}\right\}\left[g(Y,Z)\xi-\eta(Z)Y\right],
\end{equation}
\begin{equation} \label{GrindEQ__3_3_} 
\widetilde{C}(X,Y)\xi=\left\{\kappa-\frac{r}{2n(2n+1)}\right\}\left[\eta(Y)X-\eta(X)Y\right].
\end{equation}
\noindent The condition $\widetilde{C}(\xi,X)\cdot\widetilde{C}=0$ implies that
\begin{equation} \label{GrindEQ__3_4_} 
\widetilde{C}(\xi,U)\widetilde{C}(X,Y)\xi-\widetilde{C}(\widetilde{C}(\xi,U)X,Y)\xi-\widetilde{C}(X,\widetilde{C}(\xi,U)Y)\xi=0.
\end{equation}
\noindent Using (3.1) and  (3.2) in (3.3), we have
\begin{equation} \label{GrindEQ__3_5_} 
\left\{\kappa-\frac{r}{2n(2n+1)}\right\}\left[\widetilde C(X,Y)U-\left(\kappa-\frac{r}{2n(2n+1)}\right)(g(Y,U)X-g(X,U)Y)\right]=0.
\end{equation}
\noindent So either $r=2n\kappa(2n+1)$, or
\begin{equation} \label{GrindEQ__3_6_} 
\widetilde C(X,Y)U=\left(\kappa-\frac{r}{2n(2n+1)}\right)(g(Y,U)X-g(X,U)Y)=0.
\end{equation}
If $r=2n\kappa(2n+1)$, comparing the value of $r=2n(2n-2+\kappa)$, we get $\kappa=(1-\frac{1}{n})$. Hence, the manifold is $N(1-\frac{1}{n})$-contact metric manifold. That is, it is  locally isometric to Example 2.1.

\noindent Also from (1.5) and (3.5), we obtain
\begin{equation} \label{GrindEQ__3_7_} 
R(X,Y)U=\kappa \left\{g(X,U)Y-g(Y,U)X\right\}.
\end{equation}
This implies that $M$ is of constant curvature $(-\kappa)$. Consequently it is locally isometric to the hyperbolic space $H^{2n+1}(-\kappa)$. Conversely, if the manifold is $H^{2n+1}(-\kappa)$, then (3.6) holds, which yields $r=2n\kappa(2n+1)$. Then from (3.1) it follows that $\widetilde C(\xi,X)=0$. Again if $k=(1-\frac {1}{n})$, then in view of (2.14) we have $r=2(n-1)(2n+1)$ and hence it follows from (3.1) that $\widetilde C(\xi,X)=0$. Consequently, $\widetilde C(\xi,X)\cdot R=0$ since $\widetilde C(\xi,X)$ acts as a derivation. Hence $\widetilde C(\xi,X)\cdot \widetilde C=0$. This proves the theorem.

\begin{corollary}
\label{cor5.8}
A $(2n+1)$-dimensional $N(\kappa)$-contact metric manifold $M$ $(n>1)$, satisfies
\begin{equation*} \label{GrindEQ__3_8_} 
\widetilde{C}(\xi,X)\cdot R=0
\end{equation*}
if and only if either the manifold is $N(1-\frac{1}{n})$-contact metric manifold or it is locally isometric to the hyperbolic space $H^{2n+1}(-\kappa)$.
\end{corollary}
\begin{corollary}
\label{cor5.7}
In a $(2n+1)$-dimensional $N(\kappa)$-contact metric manifold $M$ $(n>1)$, we have $R\cdot \widetilde{C}=R\cdot R$.
\end{corollary}
\noindent In particular, if we consider a $3$-dimensional $N(\kappa)$-contact metric manifold, then $n=1$ and in that case, we have $\kappa=0$. Hence, in view of Proposition 2.2, we have the following result.
\begin{corollary}
\label{cor5.9}
A $3$-dimensional $N(\kappa)$-contact metric manifolds satisfies $\widetilde{C}(\xi,X)\cdot\widetilde{C}=0$ if and only if the manifold is flat.
\end{corollary}
\noindent In\cite{BLA5} it is proved that a $3$-dimensional $N(\kappa)$-contact metric manifolds is either Sasakian, flat or locally isometric to a left invariant metric on the Lie group $SU(2)$ or $SL(2,\Re)$. Hence we have the following result.

\begin{corollary}
\label{cor5.9}
If a $3$-dimensional $N(\kappa)$-contact metric manifold satisfies the condition $\widetilde{C}(\xi,X)\cdot\widetilde{C}=0$ then the manifold is either Sasakian, flat or locally isometric to a left invariant metric on the Lie group $SU(2)$ or $SL(2,\Re)$.
\end{corollary}

\begin{corollary}
\label{cor5.9}
If a $3$-dimensional $N(\kappa)$-contact metric manifold satisfies the condition $\widetilde{C}(\xi,X)\cdot R=0$ then the manifold is either Sasakian, flat or locally isometric to a left invariant metric on the Lie group $SU(2)$ or $SL(2,\Re)$.
\end{corollary}

\begin{theorem}
\label{thm4.11}
A $(2n+1)$-dimensional $N(\kappa)$-contact metric manifold $M$ $(n>1)$, satisfies
\begin{equation*} \label{GrindEQ__3_9_} 
\widetilde{C}(\xi,X)\cdot S=0
\end{equation*}
if and only if either the manifold is $N(1-\frac{1}{n})$-contact metric manifold or it is an Einstein manifold.
\end{theorem}

\noindent\textbf{Proof.} The condition $\widetilde{C}(\xi,X)\cdot S=0$ implies that
\begin{equation} \label{GrindEQ__3_10_} 
S(\widetilde{C}(\xi, X,Y),\xi)+S(Y,\widetilde{C}(\xi,X)\xi)=0.
\end{equation}
\noindent With the help of (3.1), we get from (3.7) that
\begin{equation} \label{GrindEQ__3_11_}
\left\{\kappa-\frac{r}{2n(2n+1)}\right\} \left[g(X,Y)S(\xi,\xi)-S(X,\xi)\eta(Y)+S(Y,\xi)\eta(X)-S(X,Y)\right]=0.
\end{equation}
Making use of (2.10), we have
\begin{equation} \label{GrindEQ__3_12_}
\left\{\kappa-\frac{r}{2n(2n+1)}\right\} \left[S(X,Y)-2nkg(X,Y)\right]=0.
\end{equation}
From (3.9), we obtain that either $r=2n\kappa (2n+1)$ or, $S(X,Y)=2n\kappa g(X,Y)$. If $r=2n\kappa (2n+1)$, comparing the value of $r=2n(2n-2+\kappa)$ we get $\kappa=(1-\frac{1}{n})$. Hence, the manifold is $N(1-\frac{1}{n})$-contact metric manifold. This completes the proof.\\

\noindent In particular, if we consider a $3$-dimensional $N(\kappa)$-contact metric manifold, then $n=1$ and in that case, we have $\kappa=0$. Hence, in view of Proposition 2.2, we have the following result.
\begin{corollary}
\label{cor2.1}
A $3$-dimensional $N(\kappa)$-contact metric manifolds satisfies $\widetilde{C}(\xi,X)\cdot S=0$ if and only if the manifold is flat.
\end{corollary}
Therefore, a $N(\kappa)$-contact metric manifold satisfying $\widetilde{C}(\xi,X)\cdot S=0$ is an Einstein manifold. Therefore, in view of Proposition 2.4, the manifold is a Sasakian manifold. In view of the above discussions, we can state the following:

\begin{theorem}
\label{thm4.11}
A $(2n+1)$-dimensional $N(\kappa)$-contact metric manifold $M$ $(n>1)$, satisfies
\begin{equation*} \label{GrindEQ__3_9_} 
\widetilde{C}(\xi,X)\cdot S=0
\end{equation*}
if and only if the manifold is an Einstein-Sasakian manifold.
\end{theorem}
\noindent Next, we have following theorem

\begin{theorem}
\label{thm4.2}
A $(2n+1)$-dimensional $N(\kappa)$-contact metric manifold $M$ $(n>1)$, satisfies
\begin{equation*} \label{GrindEQ__3_13_} 
\widetilde{C}(\xi,X)\cdot C=0
\end{equation*}
if and only if either the manifold is $N(1-\frac{1}{n})$-contact metric manifold or it is an $\eta$-Einstein manifold.
\end{theorem}

\noindent\textbf{Proof.} The condition $\widetilde{C}(\xi,X)\cdot C=0$ implies that

\begin{equation} \label{GrindEQ__3_14_} 
\widetilde{C}(\xi, U)C(X,Y)Z-C(\widetilde{C}(\xi,U)X,Y)Z-C(X,\widetilde{C}(\xi,U)Y)Z=0.
\end{equation}
\noindent By virtue of (3.1), we get from (3.10) that
\begin{equation} \label{GrindEQ__3_15_}
\left\{\kappa-\frac{r}{2n(2n+1)}\right\} \left[g(C(X,Y)Z,U)\xi-\eta(C(X,Y)Z)U-g(U,X)C(\xi,Y)Z\right.
\end{equation}
$$\,\,\,\,\,\,\,\,\,\,\,\,\,\,\,\,\,\,\,\,\,\,\,\,+\eta(X)C(Y,U)Z-g(U,Y)C(X,\xi)Z+\eta(Y)C(X,U)Z$$
$$\left.{\,\,\,\,\,\,\,\,\,\,\,\,\,\,\,\,\,\,\,\,\,\,\,\,\,}+\eta(Z)C(X,U)Y-g(U,Z)C(X,Y)\xi\right]=0.$$
\noindent Thus either $r=2n\kappa(2n+1)$ or, 
\begin{equation} \label{GrindEQ__3_16_}
\left[g(C(X,Y)Z,U)\xi-\eta(C(X,Y)Z)U-g(U,X)C(\xi,Y)Z\right.
\end{equation}
$$+\eta(X)C(Y,U)Z-g(U,Y)C(X,\xi)Z+\eta(Y)C(X,U)Z$$
$$\left.+\eta(Z)C(X,U)Y-g(U,Z)C(X,Y)\xi\right]=0.$$

\noindent If $r=2n\kappa(2n+1)$, comparing the value of $r=2n(2n-2+\kappa)$, we get $\kappa=(1-\frac{1}{n})$. Hence, the manifold is $N(1-\frac{1}{n})$-contact metric manifold. Again, taking the inner product of (3.12) with $\xi$, we obtain
\begin{equation} \label{GrindEQ__3_17_}
\left[C(X,Y,Z,U)-\eta(C(X,Y)Z)\eta(U)-g(U,X)\eta(C(\xi,Y)Z)\right.
\end{equation}
$$+\eta(X)\eta(C(Y,U)Z)-g(U,Y)\eta(C(X,\xi)Z)+\eta(Y)\eta(C(X,U)Z)$$
$$\left.+\eta(Z)\eta(C(X,U)Y)-g(U,Z)\eta(C(X,Y)\xi)\right]=0.$$

\noindent Using (1.2), (1.6), (2.6) and (2.10) in (3.13), we have
\begin{equation} \label{GrindEQ__3_18_} 
S(Y,Z)=\alpha g(Y,Z)+\beta \eta(Y)\eta(Z),
\end{equation}
where $\alpha=(\kappa-\frac{r}{2n(2n+1)})$ and $\beta=(\kappa+\frac{r}{2n(2n+1)})$.
\noindent This completes the proof.

\noindent In particular, if we consider $3$-dimensional $N(\kappa)$-contact metric manifold, then $n=1$ and in that case, we have $\kappa=0$. Hence, in view of Proposition 2.2, we have the following result.

\begin{corollary}
\label{cor5.10}
A $3$-dimensional $N(\kappa)$-contact metric manifold is $\widetilde{C}(\xi,X)\cdot C=0$ if and only if the manifold is flat.
\end{corollary}

\noindent Using the Proposition 2.3, we have the following:
\begin{theorem}
\label{thm4.5}
A $(2n+1)$-dimensional $N(\kappa)$-contact metric manifold satisfies the condition $\widetilde{C}(\xi,X)\cdot C=0$, $(n>1)$, then the manifold is a Sasakian manifold.
\end{theorem}
\begin{theorem}
\label{thm4.3}
If a $(2n+1)$-dimensional $N(\kappa)$-contact metric manifold $M$ $(n>1)$, satisfies the condition $C\cdot S=0$, then
\begin{equation*}\label{GrindEQ__3_19_} 
 {S^{2}(Y,Z)=\left\{(2n-1)\kappa-\frac{(2n-1)r}{2n(2n+1)}\right\}S(Y,Z)-2n\kappa \left\{\kappa+\frac{(2n-1)r}{2n(2n+1)}\right\}g(Y,Z)}.
\end{equation*}
\end{theorem}

\noindent\textbf{Proof.} In view of (1.10) the condition $C\cdot S=0$ implies that
\begin{equation}\label{GrindEQ__3_20_} 
S(C(X,Y)Z,W)+S(Z,C(X,Y)W)=0,
\end{equation}
where $X,Y,Z,W\in \chi(M)$. Taking $X=\xi$ in (3.15), we get
\begin{equation}\label{GrindEQ__3_21_} 
S(C(\xi,Y)Z,W)+S(Z,C(\xi,Y)W)=0.
\end{equation}
By virtue of (1.2),(1.6),(2.6) and (2.10), equation (3.16) reduces
\begin{equation}\label{GrindEQ__3_22_} 
\left\{\kappa+\frac{r}{2n(2n+1)}\right\}\left[2n\kappa g(Y,Z)\eta(W)-S(Y,W)\eta(Z)\right.
\end{equation}
$$\left.+2n\kappa g(Y,W)\eta(Z)-S(Y,Z)\eta(W)\right]$$
$$+\frac{1}{2n-1}\left[g(Y,W)S^{2}(\xi,Z)+S^{2}(W,Z)\eta(Z)\right.$$
$$\left.\,\,\,\,\,\,\,\,\,\,\,\,\,-g(Y,Z)S^2(\xi,W)+\eta(W)S^2(Y,Z)\right]=0.$$

\noindent So, replacing $W$ with $\xi$ in (3.17) and using (2.10), we have
\begin{equation}\label{GrindEQ__3_23_} 
S^{2}(Y,Z)=\left\{(2n-1)\kappa-\frac{(2n-1)r}{2n(2n+1)}\right\}S(Y,Z)
\end{equation}
$$-2n\kappa \left\{\kappa+\frac{(2n-1)r}{2n(2n+1)}\right\}g(Y,Z).$$
\noindent This completes the proof.
\begin{lemma}\label{lemm5.1}\cite{DES4}
Let $A$ be a symmetric $(0,2)$-tensor at point $x$ of a semi-Riemannian manifold $(M^{n},g)$, $(n>1)$, and let $T=g\wedge A$ be the Kulkarni-Nomizu product of $g$ and $A$. Then the relation
\begin{equation*} \label{GrindEQ__3_24_} 
\ T\cdot T=\alpha Q(g,T),\,\,\, \alpha \in \Re 
\end{equation*}
is satisfied at $x$ if and only if the condition
\begin{equation*} \label{GrindEQ__3_25_} 
\ A^{2}=A \alpha +\lambda g, \,\,\, \lambda \in \Re 
\end{equation*}
holds at $x$.
\end{lemma}

\noindent With the help of Theorem 3.6 and Lemma 3.1 we have the following result:

\begin{corollary}
\label{cor4.1}
Let $M$ be a (2n+1)-dimensional $N(\kappa)$ contact metric manifold $M$ $(n>1)$, satisfying the condition $C\cdot S=0$. Then $T\cdot T=\alpha Q(g,T)$, where $T=g\wedge S $ and $\alpha=\left((2n-1)\kappa-\frac{(2n-1)r}{2n(2n+1)}\right)$.
\end{corollary}

\begin{corollary}
\label{cor5.11}
Let $M$ be an $(2n+1)$-dimensional $\eta$-Einstein $N(\kappa)$-contact metric manifold $M$ $(n>1)$. Then the condition $C\cdot S=0$ holds on $M$.
\end{corollary}
\noindent\textbf{Proof.} We suppose that $M (n>1)$, be an $\eta$-Einstein $N(\kappa)$-contact metric manifold. It is well-known that Weyl tensor $C$ has all symmetries of a curvature tensor. In view of (1.10) and (2.18) we have
\begin{equation}\label{GrindEQ__3_26_} 
(C(X,Y)\cdot S)(Z,W)=\left(\frac{r}{2n}-(2n+1)\kappa\right)[\eta(C(X,Y)Z)\eta(W)+\eta(C(X,Y)W)\eta(Y)],
\end{equation}
for all vector fields $X,Y,Z,W\in \chi (M)$. By using (1.6), (2.10), (2.12) and (2.18), by a straightforward calculation, we get $(C(X,Y)\cdot S)(Z,W)=0$, so we get required result.
\noindent This completes the proof.

\begin{theorem}
\label{thm4.5}
If a $(2n+1)$-dimensional $N(\kappa)$-contact metric manifold $M$ is Weyl-pseudosymmetric then $M$ is either locally isometric to the Riemannian product $E^{n+1}(0)\times S^{n}(4)$ for $n>1$ and flat for $n=1$ or $\eta$-Einstein manifold.
\end{theorem}
\noindent\textbf{Proof.} Let $M$ $(n>1)$, be a Weyl-pseudosymmetric $N(\kappa)$-contact metric manifold. Then from (1.9) we have
\begin{equation} \label{GrindEQ__3_27_} 
(R(X,Y)\cdot C(Z,U,V)=f_{C}Q(g,C)(Z,U,V;X,Y).
\end{equation}
In view of (1.7) and (1.8), equation (3.20) yields
\begin{equation}\label{GrindEQ__3_28_} 
 R(X,Y)C(Z,U,V)-C(R(X,Y)Z,U)V-C(Z,R(X,Y)U)V
\end{equation} 
$$ -C(Z,U)R(X,Y)V=f_C[(X\wedge Y)C(Z,U)V-C((X\wedge Y)Z,U)V$$
$$-C(Z,(X\wedge Y)U)V-C(Z,U)(X\wedge Y)V].$$

\noindent Taking $X=\xi$ in (3.21), using (1.2) and (2.6) we have
\begin{equation}\label{GrindEQ__3_29_} 
 \kappa[g(Y,C(Z,U)V)\xi-\eta(C(Z,U)V)]-\kappa[g(Y,Z)C(\xi,U)V-\eta(Z)C(Y,U)V]
 \end{equation}
$$ -\kappa [g(Y,U)C(Z,\xi,V)-\eta(U)C(Z,Y,V)]$$
$$-\kappa[g(Y,V)C(Z,U)\xi-\eta(V)C(Z,U,Y)]=f_C[g(Y,C(Z,U)V)\xi-\eta(C(Z,U)V)Y$$
$$-g(Y,Z)C(\xi,U)V+\eta(Z)C(Y,U)V-g(Y,U)C(Z,\xi)V$$
$$+\eta(U)C(Z,U)V-g(Y,V)C(Z,U)\xi+\eta(V)C(Z,U)Y].$$
\noindent Taking the inner product of (3.22) with $\xi$ and then putting $Y=Z$, we get
\begin{equation}\label{GrindEQ__3_30_} 
(\kappa+f_C)[g(Y,C(Z,U,V,Z)\xi-\eta(V)\eta(C(Z,U)Z)
\end{equation}
$$+g(Z,Z)\eta(C(\xi,U)V)-g(Z,U)\eta(C(Z,\xi)V)]=0.$$

\noindent On contracting (3.23) along $Z$, we obtain
\begin{equation}\label{GrindEQ__3_31_} 
(\kappa+f_C)\eta(C(\xi,U)V)=0.
\end{equation}
\noindent If $f_C=0$, that is, the $M$ is Weyl-semisymmetric. Then from (3.24), either $\kappa=0$, or
\begin{equation}\label{GrindEQ__3_32_}
\eta(C(\xi,U)V)=0,
\end{equation}
which gives 
\begin{equation}\label{GrindEQ__3_33_}
S(U,V)=\alpha g(U,V)+\beta \eta(U)\eta(V),
\end{equation}
where $\alpha=-\left\{\kappa+\frac{(2n-1)r}{2n(2n+1)}\right\}$, $\beta=\left\{(2n+1)\kappa+\frac{(2n-1)r}{2n(2n+1)}\right\}$.\newline
If $\kappa=0$, then by Proposition 2.2, the manifold is  locally isometric to the Riemannian product $E^{n+1}(0)\times S^{n}(4)$ for $n>1$ and flat for $n=1$.

\noindent On the other hand  if $f_C\neq0$ and $\eta(C(\xi,U)V)\neq0$ then (3.24) gives $f_C=-\kappa$. So we have the following result.

\begin{corollary}
\label{cor5.12}
Every Weyl-pseudosymmetric $(2n+1)$-dimensional $N(\kappa)$-contact metric manifold, $(n>1)$, is of the form $R\cdot C=-\kappa Q(g,C)$.
\end{corollary}
\begin{corollary}

\label{cor5.13}
Every Weyl-pseudosymmetric $(2n+1)$-dimensional Sasakian manifold, $(n>1)$, is of the form $R\cdot C=-Q(g,C)$.
\end{corollary}

\noindent Using the Proposition 2.3, we have the following
\begin{theorem}
\label{thm4.5}
A $(2n+1)$-dimensional Weyl-pseudosymmetric $N(\kappa)$-contact metric manifold, $(n>1)$, is a Sasakian manifold.
\end{theorem}

\section{Example}

\noindent\textbf{Example.4.1} We consider a $3$-dimensional manifold $M=\{ (x,y,z)\in \Re ^{3}, (x,y,z) \ne (0, 0, 0)\},$ where $(x,y,z)$ is the standard coordinate in $\Re ^{3} $. Let $(e_{1} , e_{2} , e_{3} )$ be linearly independent vector fields in $\Re^{3}$ defined by
\begin{equation*} \label{GrindEQ__3_34_} 
e_{1}=\frac{\partial}{\partial y},\,\,\,\,\,e_{2}=\frac{\partial}{\partial x}-2z \frac{\partial}{\partial y}+2y\frac{\partial}{\partial z},\,\,\,\,\,\,\,e_{3}=\frac{\partial}{\partial z},
\end{equation*} 
and
\begin{equation*} 
[e_{1} ,e_{2} ]=2e_3,  \,\,\,\,\,\,\, [e_{3} ,e_{1} ]=0, \,\,\,\,\,\,\,   [e_{2} ,e_{3} ]=2e_{1}.
\end{equation*} 
Let $g$ be the Riemannian metric defined by
\begin{equation*} \label{GrindEQ__3_35_} 
g(e_{1} ,e_{2} )=g(e_{2} ,e_{3} )=g(e_{1} ,e_{3} )=0, \, g(e_{1} ,e_{1} )=g(e_{2} ,e_{2} )=g(e_{3} ,e_{3} )=1.
\end{equation*} 
  
\noindent Let $\eta$ be the $1$-form such that
\[\eta (U)=g(U,e_{3} )\] 
for any $U \in \chi(M)$. Let $\phi$ be the $(1,1)$-tensor field defined by 
\[\phi e_{1} =e_{2} ,\, \, \, \, \, \phi e_{2} =-e_{1} ,\, \, \, \, \, \, \, \phi e_{3} =0.\] 
Making use of the linearity of $\phi $ and $g$, we have
\[\eta (e_{3} )=1,\] 
\[\phi ^{2} (U)=-U+\eta (U)e_{3} \] 
and
\[g(\phi U,\phi V)=g(U,V)-\eta (U)\eta (V),\] 
for any $U, W \in \chi(M).$ Moreover
\begin{equation*} \label{GrindEQ__3_36_} 
he_1=-e_1,\,\,\,\,\,he_2=2e_2,\,\,\,\,\,\,\,he_3=0.
\end{equation*}
The Riemannian connection $\nabla $ of metric tensor $g$ is given by Koszul's formula as
\begin{equation*} \label{GrindEQ__3_37_} 
\begin{array}{l} {2g(\nabla _{U} V,W)=U(g(V,W))+V(g(W,X))-W(g(U,V))} \\ {\, \, \, \, \, \, \, \, \, \, \, \, \, \, \, \, \, \, \, \, \, \, \, \, \, \,\,\,\,\,\,\,\,\,\,\, -g(U,[V,W])-g(V,[U,W])+g(W,[U,V]).} \end{array} 
\end{equation*} 
Using the Koszul's formula, we get 
\begin{equation*} \label{GrindEQ__3_38_} 
\left\{\begin{array}{l} {\nabla _{e_{2} } e_{3} =2e_{1},\, \, \, \,\,\,\,\,\, \, \, \nabla _{e_{2} } e_{2} =0,\, \, \, \,\,\,\,\,\, \,
 \,\,\,\,  \nabla _{e_{2} } e_{1} =-2e_{3},} \\ {\nabla _{e_{3} } e_{3} =0,\, \, \,\,\,\,\,\,\,\,\, \, \, \, \, \, \, \, \, \, \, \, 
\nabla _{e_{3} } e_{2} =0,\, \, \, \, \, \, \, \, \, \, \,\,\,\, \nabla _{e_{3} } e_{1} =0, } \\ {\nabla _{e_{1} } e_{3} =0,\,\,\,\,\, \, \, \, \, \, \, \, \,\,\,\,  \nabla _{e_{1} } e_{2} =0,\, \, \, \,\,\,\,\,\,\nabla _{e_{1} } e_{1} =0.} \end{array}\right.  
\end{equation*} 
Consequently, the manifold satisfies the relation
\begin{equation*}\label{GrindEQ__3_39_} 
\nabla_{e_1}e_3=-\phi e_1-\phi (he_1),
\end{equation*}
and
\begin{equation*}\label{GrindEQ__3_40_} 
\nabla_{e_2}e_3=-\phi e_2-\phi (he_2).
\end{equation*}
Thus we have
\begin{equation*}\label{GrindEQ__3_41_} 
\nabla_{X}\xi=-\phi X-\phi hX.
\end{equation*}
for any vector field $X$. Therefore, the manifold is a contact metric manifold for $e_3=\xi$.\newline Now, we find the components of curvature tensor as 
\begin{equation*} \label{GrindEQ__11_7_} 
R(e_{2} ,e_{3} )e_{3} =0 , \,\,\, R(e_{2} ,e_{1} )e_{3} =0, \,\,\, R(e_{1} ,e_{3} )e_{3} =0.
\end{equation*} 

\noindent From the expression of above, we conclude the manifold is $N(0)$-contact metric manifold.

\noindent\textbf{Example.4.2} We consider a $5$-dimensional differentiable manifold 
$$M=\{(x, y, z,u,v) \in \mathbb{R}^5 \mid (x, y, z,u,v) \ne (0, 0, 0)\},$$
 where $(x, y, z,u,v)$ denote the standard coordinate  in $\mathbb{R}^5$. Let $(e_{1},e_{2},e_{3},e_{4},e_{5})$ are five
vector fields in $\mathbb{R}^{5}$ which satisfies
\begin{equation*} 
[e_{1},e_{2}]=-\lambda e_{2},\,\,\,[e_{1},e_{3}]=-\lambda e_{3},\,\,\,\,[e_{1},e_{4}]=0,\,\,\,\,[e_{1},e_{5}]=0,
\end{equation*}
\begin{equation*} 
[e_{i},e_{j}]=0,\,\,\,where\,\,i,j=2,3,4,5.
\end{equation*}
 \noindent We also define the Riemannian metric $g$ by
\begin{equation*} 
g(e_{1},e_{1})=g(e_{2},e_{2})=g(e_{3},e_{3})=g(e_{4},e_{4})=g(e_{5},e_{5})=1.
\end{equation*}
\begin{equation*} 
g(e_{1},e_{i})=g(e_{i},e_{j})=0,\,\,for\,\,\,i\neq j;i,j=2,3,4,5.
\end{equation*}

\noindent Let the $1-$form $\eta$ be $\eta(Z)=g(Z, e_1)$ for any $Z \in \chi(M)$.\newline Let $\phi$ be the $(1,1)$-tensor field defined by
\begin{equation*}
\phi(e_1)=0,\,\,\phi(e_2)=e_4,\,\,\,\phi(e_3)=e_5,\,\,\,\phi(e_4)=-e_2,\,\,\,\phi(e_5)=-e_3.
\end{equation*}
\noindent By the linearity properties of $\phi$ and $g$, we have
$$\phi^{2}X=-X+\eta(X)e_1,  \eta(e_1)=1,  g(\phi{X}, \phi{Y})=g(X, Y)-\eta(X)\eta(Y)$$
for arbitrary vector fields $X, Y \in \chi(M)$. Moreover,
\begin{equation*}
 he_1=0, \,\, he_2=\frac{\lambda}{2}e_4,\,\,he_3=\frac{\lambda}{2}e_5,\,\,he_4=\frac{\lambda}{2}e_2,\,\,\,he_5=\frac{\lambda}{2}e_3.
\end{equation*}
We recall the Koszul's formula as
\begin{eqnarray*}
2g(\nabla_{X}Y, Z)&=&Xg(Y, Z)+Yg(X, Z)-Zg(X, Y)\\&&
-g(X, [Y, Z])-g(Y, [X, Z])+g(Z, [X, Y])
\end{eqnarray*}
for arbitrary vector fields $X, Y, Z \in \chi(M)$. Using Koszul's formula we get
\begin{eqnarray*}
&&\nabla_{e_1}e_{1}=0,\quad\nabla_{e_1}e_{2}=0,\quad\nabla_{e_1}e_{3}=0,\quad\nabla_{e_1}e_{4}=0,\quad\nabla_{e_1}e_{5}=e_1,\\&&
\nabla_{e_2}e_{1}=-e_4+\frac{\lambda}{2} e_2,\quad\nabla_{e_2}e_{2}=-\lambda e_1, \quad\nabla_{e_2}e_{3}=0,\quad\nabla_{e_2}e_{4}=0,\quad\nabla_{e_2}e_{5}=0, \\&&
\nabla_{e_3}e_{1}=-e_5+\frac{\lambda}{2} e_3,\quad\nabla_{e_3}e_{2}=0,\quad\nabla_{e_3}e_{3}=-\lambda e_1,\quad\nabla_{e_3}e_{4}=0,\quad\nabla_{e_3}e_{5}=0,\\&&
\nabla_{e_4}e_{1}=e_2-\frac{\lambda}{2}e_4, \quad\nabla_{e_4}e_{2}=0,\quad\nabla_{e_4}e_{3}=0,\quad\nabla_{e_4}e_{4}=0,\quad\nabla_{e_4}e_{5}=0,\\&&
\nabla_{e_5}e_{1}=e_3-\frac{\lambda}{2}e_5, \quad\nabla_{e_5}e_{2}=0,\quad\nabla_{e_5}e_{3}=0,\quad\nabla_{e_5}e_{4}=0,\quad\nabla_{e_5}e_{5}=0. 
 \end{eqnarray*}

\noindent With the help of above relation, it is notice that $\nabla_{X}\xi=-\phi X-\phi h X$ for $\xi=e_1$. Therefore, the manifold is a contact metric manifold with the contact structure $(\phi, \eta, \xi, g)$.

Now, we find the curvature tensors as follows
\begin{eqnarray*}
&&R(e_1, e_2)e_1=\lambda^{2}e_2,\quad  R(e_1, e_2)e_2=-\lambda^{2}e_1,\quad R(e_1, e_3)e_1=\lambda^{2}e_3,\quad R(e_1, e_3)e_3=-\lambda^{2}e_1,\\&&
R(e_1, e_4)e_1=0, \quad R(e_1, e_4)e_4=0,\quad R(e_1, e_5)e_1=0,\quad R(e_1, e_5)e_5=0, \\&&
R(e_2, e_3)e_2=-\lambda^{2}e_3,\quad  R(e_2, e_3)e_3=-\lambda^{2}e_2,\quad R(e_2, e_4)e_2=0,R(e_2, e_4)e_4=0,\\&&
R(e_2, e_5)e_2=0, \quad R(e_2, e_5)e_5=0, \quad R(e_3, e_4)e_3=0,\quad R(e_3, e_4)e_4=0,\\&&
R(e_3, e_5)e_3=0, \quad R(e_3, e_5)e_5=0,\quad R(e_4, e_5)e_4=0,\quad R(e_4, e_5)e_5=0. 
\end{eqnarray*}
In view of the expressions of the curvature tensors we conclude that the manifold is a $N(-\lambda^{2})$-contact metric manifold.

\end{document}